\input graphicx.tex
%\$ \includegraphics[height=1 in]{ninesix.pdf}  $
%\hfill $1.10$

\def \bu {{\hskip -.1in}}

\def \CF {{\cal F}}
\def \page {{\vfill\supereject}}

%\voffset=1truein  tells how low in the page should we start!
%\hoffset=.5truein %tells how far to right should we start!
%\parindent=.5truein
%\hfuzz=3.44182pt
%\hsize 6.5truein

\def\tttt #1{{\textstyle{#1} }}

\def \CM {{\cal M}}

\def \tx {\textstyle}

\def\tttt #1{{\textstyle{#1} }}

\def \CL {{\cal L}}

\def \page{{\vfill\supereject}}

\def \magstep#1 {\ifcase#1 1000\or 1200\or 1440\or 1728\or 2074\or 2488\fi\relax}

\def\la{{\lambda}}

\def \CD {{\cal D}}

\def \CH {{\cal H}}

\def \CF {{\cal F}}

\overfullrule=0pt
\baselineskip 14pt
%\settabs 8 \columns
%\+ \hfill&\hfill \hfill & \hfill \hfill   & %%\hfill \hfill  & \hfill \hfill   &
 %& \hfill \hfill\cr
%\voffset=1truein  tells how low in the page should we start!
%\hoffset=.5truein %tells how far to right should we start!
\parindent=.5truein
\hfuzz=3.44182pt
\hsize 6.5truein
\font\ita=cmssi10  
\font\small=cmr6
\font\title=cmbx10 scaled\magstep2
\font\normal=cmr10 
\font\small=cmr6

\font\bol=cmbx12
\def\del {\partial}

\def\sig{\sigma}

\def \-> {\rightarrow}
\def\LL{\big\langle}

\def\RR {\big\rangle}

\def\DD {\Delta}

\def\la {\lambda}
\def\La {\Lambda}
\def \RA {\rightarrow}

\def\xon {x_1,x_2,\ldots ,x_n}

\def \sas {\vskip .06truein}
\def \ssas {\vskip .03truein}
\def\sa{{\vskip .125truein}}

\def\sap{{\vskip .25truein}}

\def \ses {\,=\,}
\def \sps {\, + \,}

\def \sms {\, - \,}

\def \scs {\, , \,}
\def \ess {\enskip}

\def \ssp {\hskip .25em}
\def \bigsp {\hskip .5truein}
\def \part {\vdash}

\def \DD {\Delta}

\normal

\vsize=9.5truein
%\normalbaselineskip=12pt 18pt for space and a half,
%\normalbaselines24pt for double space
\sap
\def\today{\ifcase\month\or
January\or February\or March\or April\or may\or June\or
July\or August\or September\or October\or November\or
December\fi
\space\number\day, \number\year}

\def \RA {{ \rightarrow }}

\def \CF {{\cal F}}

\def \BQ {{\bf Q}}

\def \II {1}

\def \TH {{\widetilde H}}

\def\xon {x_1,x_2,\ldots ,x_n}
\def \yon {y_1,y_2, \ldots ,y_n}

\font\small=cmr6
\def \scs {\ssp , \ssp}
\def \ess {\enskip}
\def \ssp {\hskip .25em}
\def \bigsp {\hskip .5truein}
\def \part {\vdash}
\font\title=cmbx10 scaled\magstep2
\font\normal=cmr10 
\def\today{\ifcase\month\or
January\or February\or March\or April\or may\or June\or
July\or August\or September\or October\or November\or
December\fi
\space\number\day, \number\year}
\headline={\small  
 A. M. Garsia, M. Zabrocki \small $\hskip .36in$\hfill\hfill A Basis for the Diagonal Harmonics Alternants  \hfill \hfill \hfill\today $\ess\ess\ess$    \folio }
 \footline{\null}

\def \ux {\underline x}
\def \uy {\underline y}

\def \II {{\rm I}}

\def \CL {{\cal L}}

\def \and {\ess\ess\ess\hbox{and}\ess\ess}

\hsize 6.5 in
\vsize 9 in

\centerline{\bol 
A basis for 
the Diagonal Harmonic Alternants  }

\centerline{\bf by}

\centerline{\bf A. M. Garsia\footnote{$(\dag)$}{$  Garsia\ess 
  \ess is\ess supported\ess  by\ess an \ess
   NSF   \ess Grant$} and M. Zabrocki}
\sas

\hsize 6.5in
\noindent{\bf Abstract}

It will be shown here that there are differential operators $E,F$ and $H=[E,F]$ for each $n\ge 1$, acting on Diagonal Harmonics, yielding that $DH_n$ is a representation of $sl[2]$  (see [3] Chapter 3). 
Our main effort here is to use $sl[2]$ theory  to predict a basis for the Diagonal Harmonic Alternants, $DHA_n$. 
It can be shown that the irreducible  representations $sl[2]$ are all of the form 
$P,EP,E^2P,\cdots,E^kP$, with $FP=0$ and $E^{k+1}P=0$. The polynomial $P$ is known to be called a ``String Starter''. From  $sl[2]$
theory it follows that $DHA_n$ is a direct sum of strings.  Our main result so far is a formula for the number of string starters. A recent paper by Carlsson and Oblomkov (see [2]) constructs a basis for the space of Diagonal Coinvariants by Algebraic Geometrical tools. It would be interesting to see if any our results can be derived from theirs.
\sa

\noindent{\bf Introduction}

We set $X_n=\xon$ and $Y_n=\yon$, we will be working here with polynomials $P(X_n;Y_n)$ with rational coefficients, that is
 $P(X_n;Y_n)\in\BQ[\xon;\yon]$.

 The diagonal action of $S_n$ is defined by  setting for any $\sig \in S_n$
$$
\sig\, P(X_n;Y_n)\ses  
P(x_{\sig_1},x_{\sig_2},\ldots ,x_{\sig_n};
y_{\sig_1},y_{\sig_2},\ldots ,y_{\sig_n})\, . 
\eqno \II.1
$$
Another important tool in studying $S_n$ modules
that are  invariant under the diagonal action
 is the scalar product
$$
\LL P\scs Q \RR\ses 
L_o P(\del X_n ;\del Y_n )Q ( X_n ;Y_n).
\eqno \II.2 
$$ 
where the differential operator $P(\del X_n ;\del Y_n )$ is obtained by the replacements 
$x_i\RA \del_{x_i}$ and $y_i\RA \del_{y_i}$.
It is easy to see that we have
$$
\LL\sig  P\scs \sig Q \RR\ses \LL P\scs Q \RR
\ess\ess\ess\ess\ess\ess\ess
\hbox{(for all $\sig \in S_n$
)}.
\eqno \II.3 
$$
In this paper we will study the $S_n$ module $DH_n$ of Diagonal Harmonic polynomials. This 
module was originally   defined  as the orthogonal complement, with respect to the scalar product in I.2, of the ideal of polynomials that are invariant under the diagonal action. By a result of Hermann  Weyl (see [16])  it follows that 
$P(X_n;Y_n)\in DH_n$  if and only if
$$
 \sum_{i=1}^n\del_{x_i}^p\del_{y_i}^q\,P(X_n;Y_n)\ses 0
\bigsp \hbox{(for all $p+q\ge 1$)}\, .
\eqno \II.4 
$$
This simpler definition makes it obvious that $DH_n$ is invariant under the diagonal action.

It also immediately  follows from I.4 that if $P\in DH_n$ then all the bi-homgeneous components of $P$ are in $DH_n$.
This implies that we have the direct sum decomposition
$$
DH_n\ses \bigoplus_{0\le r+s\le {n \choose 2}}
\CH_{r,s} \big( DH_n\big),
\eqno \II.5 
$$
where $\CH_{r,s}(DH_n)$ is the  subspace of diagonal Harmonics  polynomials which are bi-homogeneous of degree $r$ in the $x's$ and    degree  $s$ in the $y's$. It then  follows that the 
 character resulting from the diagonal action of $S_n$ on $DH_n$ can be written in the form
$$
\chi^{DH_n}\ses  \sum_{0\le r+s\le {n \choose 2}} 
t^r q^s \chi^{r,s},
\ess\ess\ess\ess\ess\ess
\hbox{\big(where $\chi^{r,s}$ is the character of $\CH_{r,s}( DH_n)$
\bu\big).}
\eqno \II.6 
$$
The upper bound ${n \choose 2}$ in I.5 follows from 
the operator conjecture (proved by Mark Haiman in [13]). $DH_n$  can be obtained by applying  differential operators to the
 Vandermonde determinant $\prod_{1\le i<j\le n}(x_i-x_j)$.
\supereject

The Frobenius map ``$\CF$ ''(see [3]) considerably simplifies the operation of computing the character  
of an $S_n$ representation. Frobenius uses  the dimension  equality between the class functions of $S_n$ and the space $\La^{=n}$ of homogeneous symmetric functions of degree $n$.  This given, $\CF$ maps Class Functions onto the power basis by the formula 

\vskip -.3in
$$
\CF C_\mu \ses p_\mu/z_\mu.
\eqno \II.7
$$

\vskip -.08in
\noindent
where $C_\mu$  is the sum  of all the permutations
of cycle structure $\mu$, $p_\mu=p_1p_2\cdots p_{\l(\mu)}$, where
$\l(\mu)$ denotes the length of $\mu$ and
$$
z_\mu\ses 1^{a_1}2^{a_2}\cdots n^{a_n} 
a_1!a_2!\cdots a_n!
\ess\ess\ess\ess\ess\ess\ess\ess
\hbox{(when $ \mu = 1^{a_1}2^{a_2}\cdots n^{a_n}
\part n$)}
\eqno \II.8
$$

\vskip -.02in
\hfill$ \hbox {\includegraphics[width=1.7in]{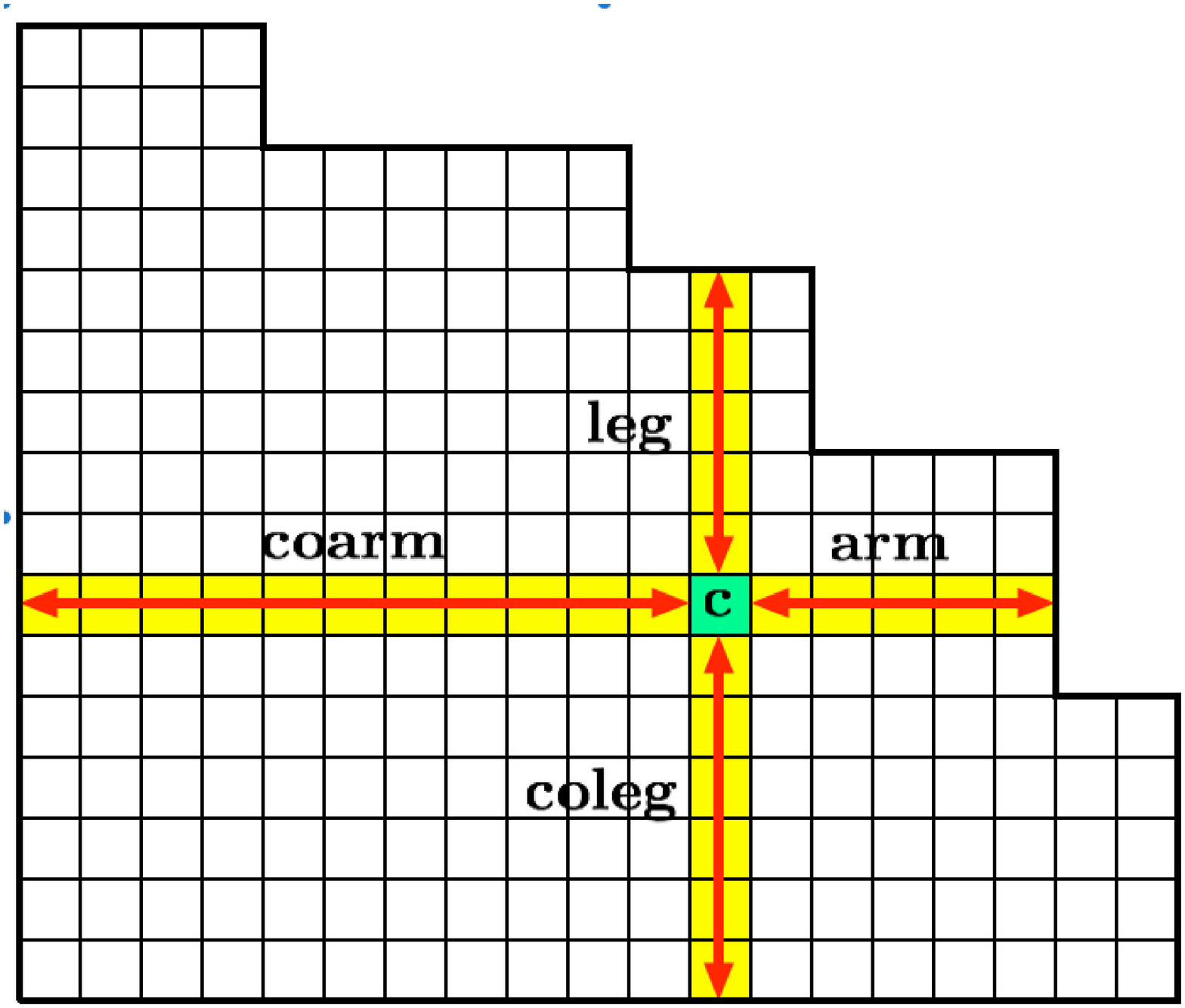}}  $

\hsize=4.6in
\vskip -1.4in

\noindent
It follows from I.7 that Young's irreducible  character $\chi^\la$
is given by the identity

\vskip -.18in
$$
\chi^\la\ses \CF^{-1} s_\la
\eqno \II.9
$$

\vskip -.08in
\noindent
If $\CM_n$ is an $S_n$  module with character 
$\chi^{\CM_n}$ then the {\ita Frobenius  Characteristic} of $\CM_n$ is the symmetric polynomial
$
\CF\, \chi^{\CM_n}
$.
It was conjectured  (in [6]) and shown by Mark Haiman using tools of Algebraic Geometry (in [12]) that the Frobenius characteristic of $DH_n$ is the symmetric rational function

\hsize=6.5in
\noindent

$$
\CF\,\chi^{DH_n}\ses \sum_{\mu\part n} {T_\mu \TH_\mu(X;q,t)
M B_\mu(q,t) \Pi_\mu(q,t)
\over w_\mu(q,t)}.
\eqno \II.10
$$
To define the ingredients that appear in this  formula   it is convenient to use Macdonald's notation [14]
 yet identify partitions with their French Ferrers diagram. Given
a partition $\mu$ and a cell $c\in \mu$, as indicated in the above display, we will introduce  four
parameters
$l_\mu(c)$, $l'_\mu(c)$, $a_\mu(c)$ and  $a'_\mu(c)$ 
called 
{\ita leg, coleg, arm } and {\ita coarm } which  give the 
number of lattice cells of $\mu$  strictly  { North},  
{ South}, 
{East } and  {West } of $c$.
 Denoting by $\mu'$ the conjugate of $\mu$, 
we set
$$
\eqalign{
 n(\mu)\ses \sum_{c\in \mu}{l'_\mu(c)} ,
 \ess\ess\ess\ess\ess 
 T_\mu=t^{n(\mu)}q^{n(\mu')}
, & \ess\ess\ess\ess M=(1-t)(1-q),
\ess\ess\ess\ess
B_\mu(q,t)=\sum_{c\in \mu}t^{l'_\mu(c)}q^{a'_\mu(c)},
 \cr
 \Pi_\mu(q,t)=\hskip-.2in\prod_{c\in \mu;c\neq(0,0)}\hskip-.2in
(1-{l'_\mu(c)}q^{a'_\mu(c)}),\ess\ess\ess\ess
& w_\mu(q,t)=\prod_{c\in \mu}(q^{a _\mu(c)} -t^{l _\mu(c)+1})(t^{l _\mu(c)} -q^{a _\mu(c)+1}).
\cr
}
\eqno \II.11
$$

\vskip -.16in
\noindent
This accounts for every thing that occurs in 
I.10 except for the modified Macdonald Basis element
$\TH_\mu(X;q,t)$. This symmetric polynomial  was conjectured in [5] and proved by Mark Haiman in [12] to be the Frobenius Characteristic of the linear span of derivatives of the alternant that corresponds to the partition $\mu$.

Our  first goal is to construct operators
 that preserve $DH_n$. We will prove this property for all the  differential operators 

\vskip -.3in
$$
a)\ess\ess F_{r,s}\ses\sum_{i=1}^n \ux_i\del_{x_i}^r\del_{y_i}^s\,,
\bigsp\ess\ess b)\ess\ess E_{r,s}\ses\sum_{i=1}^n 
y_i\del_{x_i}^r\del_{y_i}^s\,,
\ess\ess\ess\ess\ess
\hbox{(for $r+s\ge 1$
)
}.
\eqno \II.12
$$

\vskip -.12in
\noindent
We will also prove the relations
$$
\eqalign{
a)\ess\ess \big[F_{p,q}\scs F_{r,s}\big]&\ses(p-r)F_{p+r-1,q+s}\,,
\cr
b)\ess\ess \big[F_{p,q}\scs E_{r,s}\big]&\ses qF_{p+r,q+s-1}\sms rE_{p+r-1,q+s}\,,
\cr
c)\ess\ess \big[E_{p,q}\scs E_{r,s}\big]&\ses (q-s)E_{p+r,q+s-1} \,.\cr
}
\eqno \II.13
$$
Furthermore by setting
$$
a)\ess\ess F\ses F_{0,1}\,, 
\bigsp 
b)\ess\ess E\ses E_{1,0}\,,
\bigsp
c)\ess\ess H\ses [E, F]\ses
\sum_{i=1}^n
\big(\uy_i\del_{y_i}-\ux_i\del_{x_i}\big)
\eqno \II.14
$$
we derive that $DH_n$ is a direct sum of irreducible representations of $sl[2]$.
\supereject

Finally, using the results of the $q,t$-Catalan paper (see [8],[9]) we derive that the Hilbert polynomial 

\vskip -.18in
$$
c_n(q,t)\ses\CF\,\chi^{DHA_n}\Big|_{s[1^n]}
\eqno \II.15
$$

\vskip -.0921in
\noindent
of the Diagonal Alternants can be obtained by a purely combinatorial construction. To describe this construction we need further definitions. By a Dyck path in the $n\times n$ lattice square $\CL_n$ we mean a path which proceeds by $n$ unit North steps and $n$ unit East steps and goes from $(0,0)$ to $(n,n)$ always remaining 
weakly

 \hfill $\includegraphics[height=1.3 in]{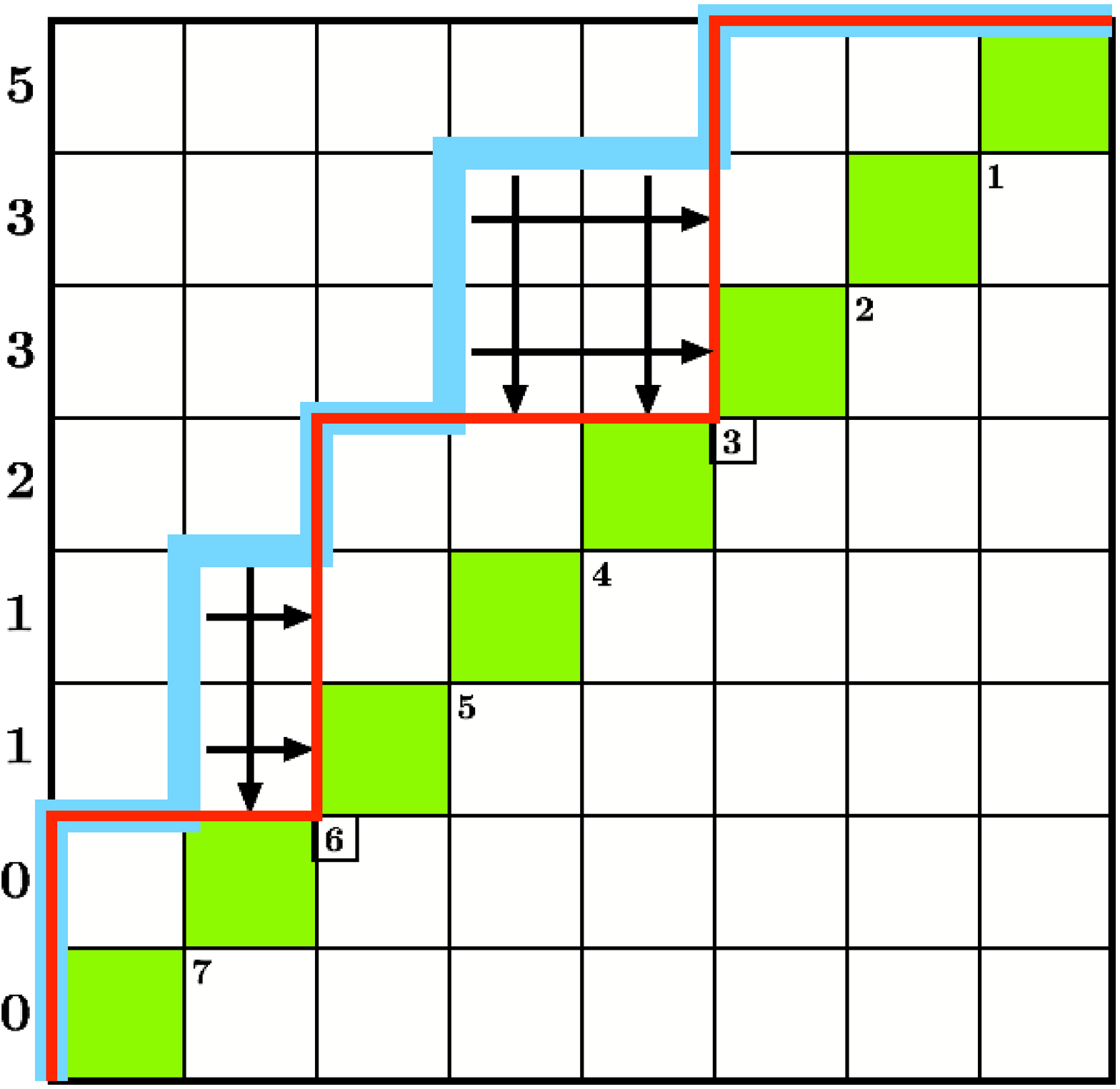}$ 

\vskip -1.36in
\hsize=5in 
\noindent
above the main diagonal of $\CL_n$. This is the straight line that joins  $(0,0)$ to $(n,n)$. In the illustration on the right we colored light green all the cells bisected by the main diagonal of $\CL_8$. A Dyck path
is depicted there in light blue. It is clear that we only need to give the abscissas of the North steps of the Dyck path $D$. Those are the integers that are on the left of the rows of $\CL_8$. Thus $D=[0,0,1,1,2,3,3,5]$.
Each $Dyck$ path has two statistics which we call $area(D)$ and $bounce(D)$. The area statistic is quite simple, its formula is $area(D)=[0,1,\ldots,n-1]-[0,d_1,\ldots,d_{n-1}]$.

\hsize=6.5in
\noindent
This is the number of cells between the Dyck path  and the lattice diagonal (the green cells). In our case $area(D)=13$.
The bounce statistic is the sum of the places where the bounce path hits the main diagonal of $\CL_n$. In our display we depicted the bounce path in red using a thinner line. In the general case it starts  straight North until it touches the West end of an East step. Then it goes straight East until it touches the diagonal. Then goes straight North until it touches the West end of an East step...  alternating straight North and straight East   until it reaches $(n,n)$. We place in the possible diagonal touching points the labels $1,2,\ldots,n-1$ as  indicated in our display.  In our example, $bounce(D)=3+6$. We must emphasize that the bounce path does not change direction by touching the East end of an East step. In our display that happens in the $5^{th}$ and $8^{th}$ rows. 

In particular we obtain 

\vskip -.19in
$$
c_n(q,t)
 \ses \sum_{D\in \CD_n}t^{bounce(D)}q^{area(D)}\ses 
 \sum_{D\in \CD_n}t^{area(D)}q^{dinv(D)}
\eqno \II.16 
$$  

\vskip -.05in
\noindent
The first identity was conjectured by Jim Haglund the second was conjectured by Mark Haiman. The {\ita dinv} statistic has also a purely combinatorial definition. 
Let $dinv_a(D)=\sum_{1\le i<j\le n}\chi(u_i=u_j)$ where $u_i$ is the contribution to the area statistic by the $i^{th}$ north step. Similarly we let $dinv_b(D)=\sum_{1\le i<j\le n}\chi(u_i=u_j+1) $, then set $dinv(D)=dinv_a(D)+dinv_b(D)$.
The problem to construct a further statistic that combined with  area gives I.16 was stated in [6], these two solutions were discovered quite a few years later.  

For $n=3$ we get

\vskip -.16in
$$
\CF DH_3
 \Big|_{s[1^3]}\ses \sum_{D\in \CD_3}t^{bounce(D)}q^{area(D)}\ses t^3\sps\ t^2 q\sps t q\sps t q^2\sps q^3 .
\eqno \II.17 
$$

\vskip -.09in
\noindent
This  identity reveals that the alternating character $\chi^{[1,1,1]}$ occurs in bi-degrees $(3,0),(2,1),(1,2),(0,3)$ and $(1,1)$. In this particular case these are Frobenius images of the $sl[2]$ strings generated by the following two alternants

\vskip -.2in 
$$ 
a)\ess\ess\DD_{1,1,1}\ses\det \pmatrix{1 & 1 & 1\cr x_1 & x_2 & x_3 \cr  x_1^2 & x_2^2 & x_3^2 \cr },
\bigsp\ess\ess 
b)\ess\ess\DD_{2,1}\ses\det \pmatrix{1 & 1 & 1\cr x_1 & x_2 & x_3 \cr y_1 & y_2 & y_3 \cr }.
\eqno \II.18 
$$
It stands to reason that  we should be able to 
construct a basis for the Diagonal Harmonics Alternants using the operators in I.12, $sl[2]$ theory and the combinatorics of the $q,t$-Catalan.
\sas

 Mark Haiman's proof in [13] of the Operator Conjecture implies
the following result
\sas

\noindent{\bf Theorem I.1}

{\ita For any $n\ge 1$ let $m$  be the dimension of $DHA_n$ in bi-degree
$(a,b)$, then for this bi-degree, 
 it is always possible to construct  $m$ sequences $1\le r_1\le r_2\le \cdots \le r_b\le n$ such that the Diagonal Harmonic Alternants
$E_{r_1,0}E_{r_2,0}\cdots E_{r_b,0}\DD_{1^n}$ are linearly independent.
 This requires that
${n \choose 2}-(r_1+r_2+\cdots + r_b)=a$
}

The only significant results we will prove in this paper, besides introducing $sl[2]$ theory and proving the invariance under the diagonal action  of the differential operators in I.12 and proving
the properties in I.13, are a formula for the number of starters and an algorithm that gives that number for every $n$.
\sas

\noindent{\bf 1. The differential operators.}

We start with an auxiliary fact concerning the interaction between multiplication operators and differential operators.
\sas

\noindent{\bf Propositiom 1.1}

{\ita For any variable $y$ and integer exponent $q\ge 1$ we have}
$$
\del_y^q \,\uy\ses q\, \del_y^{q-1}+\, \uy\,\del_y^q,
\ess\ess\ess
(\hbox{where ``$\uy$'' is the multiplication by  $y$ operator})
\eqno 1.1
$$
\noindent{\bf Proof}

Suppose that $P(y)$ is a polynomial in $y$. Then for $q=1$  we get
$$
\del_y  \,\uy\, P(y)\ses P(y)\sps \uy\,\del_y P(y).
\eqno 1.2
$$
Thus 1.1 is true for $q=1$. Proceeding by induction on $q$, suppose that 1.1 is true up to $q-1$. Then we have $$
\eqalign{
\del_y^q\, y\,P(y)&\ses \del_y\, \del_y^{q-1} y\,P(y)\ses 
\del_y\,(q-1)\del_y^{q-2}P(y)\sps \del_y\, y\, \del_y^{q-1}P(y)\ses\cr
&\ses
(q-1)\del_y^{q-1}P(y)\sps 
 \del_y^{q-1}P(y)\sps y\,\del_y^{q}P(y)
\ses q\, \del_y^{q-1}P(y)\sps y\,\del_y^{q}P(y).
\cr}
\eqno 1.3
$$
this proves 1.1.
\sas

As an example we will show that 

\noindent
{\bf  Theorem 1.1}

{\ita The  $sl[2]$ operators
$$
a)\ess\ess F\ses \sum_{i=1}^n \ux_i\del_{y_i}\,,
\bigsp\ess\ess 
b)\ess\ess E\ses \sum_{i=1}^n \uy_i\del_{x_i}\,,
\eqno 1.4
$$
preserve $DH_n[X_n;Y_n]$.}

\noindent
{\bf  Proof}

To this end we will first compute the bracket
$$
\big[\Pi_{p,q}^n\scs E\big]\ses
\sum_{i=1}^n\sum_{j=1}^n\big[
\del_{x_i}^p\del_{y_i}^q\scs \uy_j\del_{x_j}
\big]\,,
\ess\ess\ess\ess\ess\ess
\hbox{\big(where $\displaystyle \Pi_{p,q}^n=\sum_{i=1}^n\del_{x_i}^{p }\del_{y_i}^{q }$ \big).}
\eqno 1.5 
$$
Since for $j\neq i$ the differential and  multiplication operators commute, we only need to work with
$$
\big[\Pi_{p,q}^n\scs E \big]\ses
\sum_{i=1}^n\big[
\del_{x_i}^p\del_{y_i}^q\scs \uy_i\del_{x_i}
\big]\, .
\eqno 1.6
$$
Using 1.1 for $q\ge 1$ we obtain
$$
\del_{x_i}^p\del_{y_i}^q \uy_i\del_{x_i}
\ses \del_{x_i}^p \big(q\, \del_{y_i}^{q-1}+\, \uy_i\,\del_{y_i}^q\big)\del_{x_i}
\ses q\, \del_{x_i}^{p+1}\del_{y_i}^{q-1}\sps 
\uy_i\del_{x_i}^{p+1} \del_{y_i}^q.
\eqno 1.7
$$
We also have
$$
\uy_i\del_{x_i}\del_{x_i}^p\del_{y_i}^q 
\ses \uy_i\,\del_{x_i}^{p+1} \del_{y_i}^q
\eqno 1.8
$$
so 1.6 becomes 
$$
\big[\Pi_{p,q}^n\scs E \big]\ses
q\sum_{i=1}^n\del_{x_i}^{p+1}\del_{y_i}^{q-1}\sps\sum_{i=1}^n\uy_i\del_{x_i}^{p+1} \del_{y_i}^q\sms
\sum_{i=1}^n \uy_i\,\del_{x_i}^{p+1} \del_{y_i}^q\ses q\sum_{i=1}^n\del_{x_i}^{p+1}\del_{y_i}^{q-1}\,,
\eqno 1.9
$$
or  equivalently
$$
\Pi_{p,q}^n\, E \ses E\,\Pi_{p,q}^n\sps
q\,\Pi_{p+1,q-1}^n.
\eqno 1.10
$$
Thus applying $\Pi_{p,q}^n\, E$ to a diagonal harmonic polynomial $P(X_n;Y_n)$, gives
$$
\Pi_{p,q}^n\, E P(X_n;Y_n)\ses E\,\Pi_{p,q}^nP(X_n;Y_n)\sps
q\,\Pi_{p+1,q-1}^nP(X_n;Y_n) \ses 0.
\eqno 1.11
$$
The case $q=0$ is trivial.
Proving that the operator $E$ preserves diagonal harmonics. Working with $F$ we  will reach the same result using analogous steps. 
\supereject

We can use the same idea on  the operators
$$
a)\ess\ess F_{r,s}\ses\sum_{i=1}^n \ux_i\del_{x_i}^r\del_{y_i}^s\,,
\bigsp\ess\ess b)\ess\ess E_{r,s}\ses\sum_{i=1}^n 
y_i\del_{x_i}^r\del_{y_i}^s\,.
\eqno 1.12
$$ 
As before we can reduce the calculation to $j=i$ and work with 
\vskip -.15in
$$
\Pi_{p,q}^n F_{r,s}\ses \sum_{i=1}^n \del_{x_i}^p\ux_i\del_{x_i}^r\del_{y_i}^{q+s}
\ses
 \sum_{i=1}^n\big(p\, \del_{x_i}^{p+r-1}\tx
\del_{y_i}^{q+s}\sps \ux_i\del_{x_i}^{p+r}\del_{y_i}^{q+s}
\big).
\eqno 1.13
$$
\sas

\noindent
Likewise we have
\vskip -.25in
$$
F_{r,s}\Pi_{p,q}^n\ses
\sum_{i=1}^n\ux_i \del_{x_i}^{p+r}\del_{y_i}^{q+s}\,.
\eqno 1.14
$$
\vskip -.13in

\noindent
Thus
$$
\Pi_{p,q}^n\,F_{r,s} \ses F_{r,s}\Pi_{p,q}^n\sps
p\, \Pi_{p+r-1,q+s}^n
\eqno 1.15
$$
Now if $P(X_n;Y_n)$ is Diagonal Harmonic then
$$
\Pi_{p,q}^n\,F_{r,s}P(X_n;Y_n) \ses F_{r,s}\Pi_{p,q}^nP(X_n;Y_n)\sps
p\,\Pi_{p+r-1,q+s}^n P(X_n;Y_n)\ses 0\,.
\eqno 1.16
$$
Proving that $F_{r,s}$ preserves Diagonal Harmonics.
An analogous argument yields the same result for $E_{r,s}$. 
\sas

\noindent{\bf Theorem 1.2}

{\ita The following identities hold true for all $p+q\ge 1$ and $r+s\ge 1$. 
}
$$
\eqalign{
a)\ess\ess \big[F_{p,q}\scs F_{r,s}\big]&\ses(p-r)F_{p+r-1,q+s}\,,
\cr
b)\ess\ess \big[F_{p,q}\scs E_{r,s}\big]&\ses qF_{p+r,q+s-1}\sms rE_{p+r-1,q+s}\,,
\cr
c)\ess\ess \big[E_{p,q}\scs E_{r,s}\big]&\ses (q-s)E_{p+r,q+s-1} \,.\cr
}
\eqno 1.17
$$
\noindent{\bf Proof}

Reducing again to the case $j=i$ we can write
$$
\big[F_{p,q}\scs F_{r,s}\big]\ses
\sum_{i=1}^n\big[\ux_i\del_{x_i}^p \del_{y_i}^q\scs\ux_i
\del_{x_i}^r\del_{y_i}^s\big]
\ses 
\sum_{i=1}^n\ux_i(\del_{x_i}^p \ux_i)
\del_{x_i}^r\del_{y_i}^{q+s}
\sms
\sum_{i=1}^n\ux_i
(\del_{x_i}^r \ux_i)
\del_{x_i}^p\del_{y_i}^{q+s}\,.
\eqno 1.18
$$
Since we have
$$
a)\ess\ess \del_{x_i}^p\ux_i=  
p\, \del_{x_i}^{p-1}+ \ux_i\del_{x_i}^p
\bigsp\bigsp
b)\ess\ess\del_{x_i}^r \ux_i= 
r\, \del_{x_i}^{r-1}+ \ux_i\del_{x_i}^r 
\eqno 1.19
$$
The identity in 1.18 becomes, using 1.19  
$$
\big[F_{p,q}\scs F_{r,s}\big]
\ses \sum_{i=1}^n\ux_i(p\, \del_{x_i}^{p-1}+ \ux_i\del_{x_i}^p)
\del_{x_i}^r\del_{y_i}^{q+s}
\sms
\sum_{i=1}^n\ux_i
(r\, \del_{x_i}^{r-1}+ \ux_i\del_{x_i}^r)
\del_{x_i}^p\del_{y_i}^{q+s}\,.
\eqno 1.20
$$
Now this can be rearranged to
$$
\eqalign{
\big[F_{p,q}\scs F_{r,s}\big]
&\ses \sum_{i=1}^n\ux_i(p\, \del_{x_i}^{p-1})
\del_{x_i}^r\del_{y_i}^{q+s}
\sms
\sum_{i=1}^n\ux_i
(r\, \del_{x_i}^{r-1})
\del_{x_i}^p\del_{y_i}^{q+s}\sps
\sum_{i=1}^n(\ux_i^2-\ux_i^2)
\del_{x_i}^{p+r}\del_{y_i}^{q+s}\,.
\cr }
$$
From which we derive that
$$
\eqalign{
\big[F_{p,q}\scs F_{r,s}\big]
&\ses
 (p-r)\sum_{i=1}^n\ux_i \del_{x_i}^{p+r-1}\del_{y_i}^{q+s}
\ses (p-r)F_{p+r-1,q+s}\,.
\cr 
}
\eqno 1.21
$$
\vskip -.03in

\noindent
This proves a) of 1.17.

\supereject

\noindent
Next we work on b)  of 1.17.
Reducing to $j=i$ we can write 
$$
\big[F_{p,q}\scs E_{r,s}\big]
\ses
\sum_{i=1}^n\big[\ux_i\del_{x_i}^p \del_{y_i}^q ,
\uy_i
\del_{x_i}^r\del_{y_i}^s\big]
\ses 
\sum_{i=1}^n \ux_i
(\del_{y_i}^q \uy_i)
\del_{x_i}^{p+r}\del_{y_i}^s 
\sms
\sum_{i=1}^n\uy_i
(\del_{x_i}^r \ux_i)\del_{x_i}^p \del_{y_i}^{q+s}
\eqno 1.22
$$
Using 1.19  we get
$$
\eqalign{
\big[F_{p,q}\scs E_{r,s}\big]
&\ses 
\sum_{i=1}^n \ux_i
(q\, \del_{y_i}^{q-1}+ \uy_i\del_{y_i}^q)
\del_{x_i}^{p+r}\del_{y_i}^s 
\sms
\sum_{i=1}^n\uy_i
(r\, \del_{x_i}^{r-1}+ \ux_i\del_{x_i}^r )\del_{x_i}^p \del_{y_i}^{q+s}
\cr
&\ses 
q\sum_{i=1}^n \ux_i
\, 
\del_{x_i}^{p+r}\del_{y_i}^{q+s-1} 
\sms
r\sum_{i=1}^n\uy_i
 \del_{x_i}^{p+r-1} \del_{y_i}^{q+s}
 \sps \sum_{i=1}^n(\ux_i\uy_i-\uy_i\ux_i )\del_{x_i}^{p+r} \del_{y_i}^{q+s}
\cr
&\ses q\,F_{p+r,q+s-1}\sms r\,E_{p+r-1,q+s}\,.
}
\eqno 1.23
$$
This proves 1.17 b). 
The identity in 1.17 c) is proved the same  way we proved a).
\sa
 
\noindent
{\bf Remark 1.1}
 
Using b) of 1.17 we can prove the result  in c) of I.14.
 In fact, since $E=E_{1,0}$ and 
$F=F_{0,1}$, setting $p=0,q=1,r=1,s=0$ in
$$
\big[F_{p,q}\scs E_{r,s}\big] \ses 
qF_{p+r,q+s-1}\sms rE_{p+r-1,q+s}\,,
$$
we obtain
$$
\big[F_{0,1}\scs E_{1,0}\big] \ses F_{1,0}\sms E_{0,1}\,,
$$
Using 1.12 this gives that the $sl[2]$ operator in I.14
$$
H\ses \big[E\scs F\big] \ses  \sum_{i=1}^n 
y_i\del_{y_i}\sms \sum_{i=1}^n \ux_i\del_{x_i}\,,
\eqno 1.24
$$
on Diagonal Harmonics is none other than the Euler operator in the $y's$ minus the
Euler operator in the $x's$.
\sas

\noindent
{\bf Remark 1.2}

We will make multiple uses of the operators $E_{r,0}$, we must point out that that these operators commute regardless of the values of $r$. This is one of the consequences of the  identities in 1.17. In fact,  if $q=s$ in c) that will immediately force the commutativity of $E_{p,q}$  and $E_{r,s}$. The analogous  result can also be obtained for the  differential operators pairs in 1.17 a).
\sas

\noindent
{\bf Proof of Theorem I.1}

For $A=(a_1,a_2,\ldots,a_n)$ with $a_i\in Z_{\ge 0}=\{0,1,2,\ldots\}$ set
$$
E^A\ses E_{1,0}^{a_1}E_{2,0}^{a_2}\cdots E_{n,0}^{a_n}.
\eqno 1.25
$$
The proof in [13]  of the operator conjecture implies that if we set 
$\CH[\xon]=\CL_\del[\DD_{1^n}]$ (the ordinary Harmonics of $S_n$), then
$$
\sum_A E^A\CH[\xon]\ses \CH[\xon;\yon],
\eqno 1.26
$$
the space of Diagonal Harmonics. This implies that if $\{h_1, h_2,\ldots,h_{n!}\}$ is any basis for $\CH [\xon]$ then the elements $\{E^{A_1}{h_1}, E^{A_2}h_2,\ldots,E^{A_{n!}}h_{n!}\}$ span the Diagonal Harmonics, in particular they span the Diagonal Harmonic Alternants. If the basis is isotypical that is there are as many independent elements generating the irreducible representation character $\chi^\la$ as the dimension of this representation. Since the alternating element occurs with multiplicity one and is the Vandermonde in $\xon$ the  basis elements must be all of the form
$$
E^A\DD_{1^n}\ses E_{r_1,0}E_{r_2,0}\cdots E_{r_b,0}\DD_{1^n},
\eqno 1.27
$$
with $0\le r_1\le r_2\le \cdots \le r_b$. 
This proves the theorem.
\sas

\page

\noindent{\bf 2. Analysing the computer data for 
$\bf n=5,6$.}
\sas

Although our differential operators are constructing
Diagonal Harmonic alternants we will be guided  

\vskip .1in
\hfill $\includegraphics[height=2.5 in]{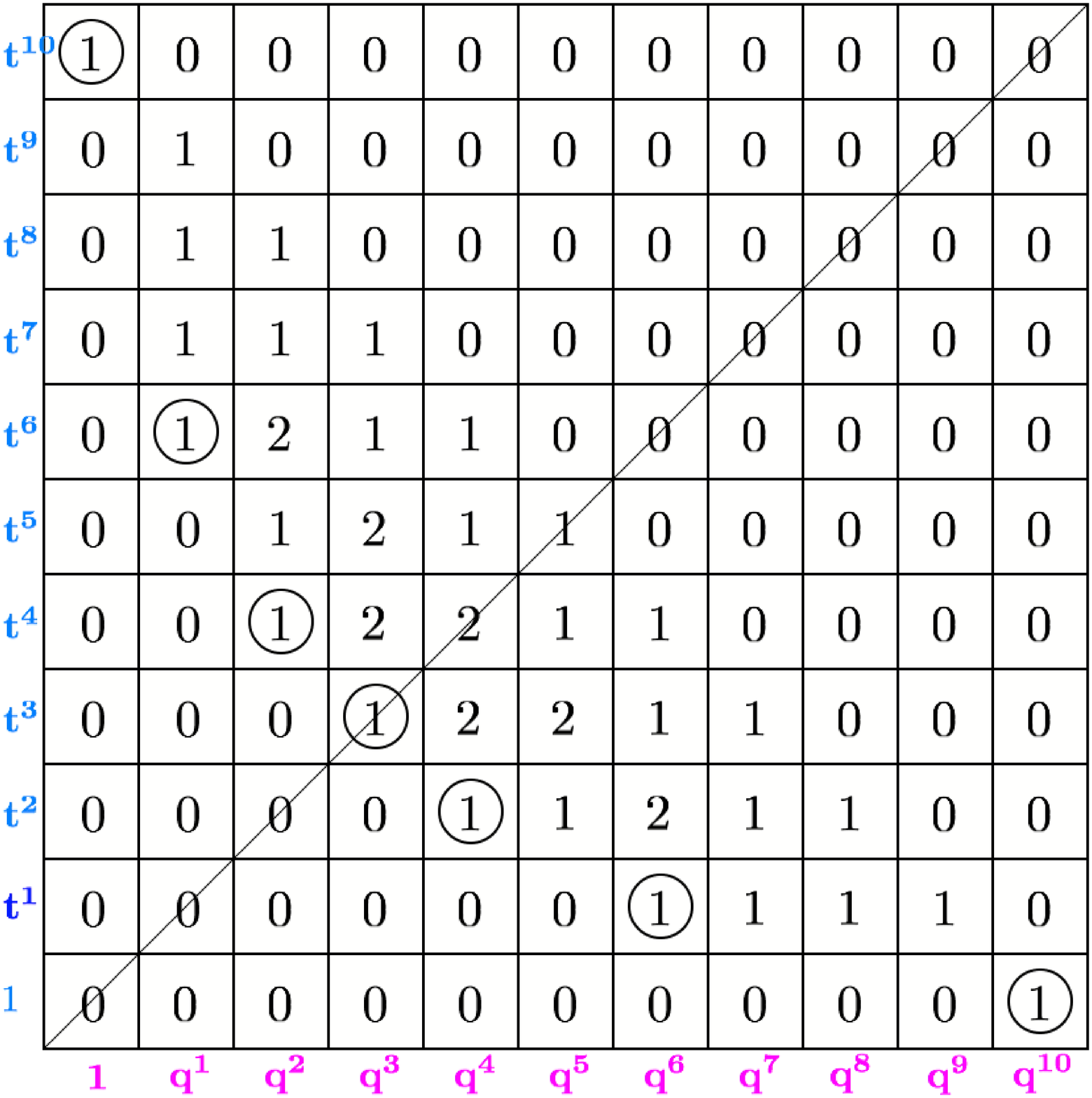}$ 

\vskip -2.64in
\hsize=3.9in
\noindent
 by the Frobenius characteristic of the alternants in $DH_n$ 
\vskip -.16in
$$
c_n(q,t)\ses \sum_{D\in\CD_n}  t^{bounce(D)}q^{area(D)}.
\eqno 2.1
$$

\vskip -.1in
\noindent
We will do that for $n=5$ and $n=6$.

In the display on the right we have depicted a visual image of this polynomial for $n=5$. For instance the term of $c_5(q,t)$ that is in the fifth row and fifth column corresponds to the Dyck path $D\in \CD_5$ for which 
$t^{bounce(D)}q^{area(D)}=2\, t^4q^4$.
We can easily locate the $sl[2]$ strings in this display. In $DH_n$ they are of the form 
$
P(X_n;Y_n)\RA EP(X_n;Y_n)\RA\cdots \RA E^kP(X_n;Y_n),
$
where $P(X_n;Y_n)$ is killed by $F$ and is a homogeneous polynomial of bi-degree $(a,b)$. 
The operator $E$ diminishes degree in $x's$ by one and increases degree in 
 $y's$ by one, thus the polynomial $E^iP(X_n;Y_n)$ is homogeneous of bi-degree $(a-i,b+i)$. The final element
is of 

\hsize 6.5 in
\noindent
bi-degree $(b,a)$ and is killed by $E$.
 The circled entries
correspond to alternants indexed by partitions. The one corresponding
to $t^{10}$ is the Vandermonde
$ \DD_{1^5}$. The one corresponding
to $t^{6}q^1$ is  $ \DD_{2 1^3}$.
$\DD_{221},\DD_{311}$ and $\DD_{32}$
form an $sl[2]$ string.
Our list ends with $\DD_{41}$ and the Vandermonde
$\DD_{5}$. This display shows that the alternants in $DH_5$ factor into 
a direct sum of
 $7$ of $sl[2]$ strings. The one of highest total degree is started by 
$\DD_{1^5}$. Let us call it {\ita String $1$}. {\ita String $2$} starts
with bi-degree $(8,1)$ and ends in bi-degree $(1,8)$. String $3$  starts
with bi-degree $(7,1)$ and ends in bi-degree $(1,7)$. String $4$  starts
with bi-degree $(6,2)$ and ends in bi-degree $(2,6)$. That accounts for the $2's$ in the image of its path. String $5$  starts
with bi-degree $(6,1)$ and ends in bi-degree $(1,6)$. These two ends are
$\DD_{2 1^3}$ and $\DD_{41}$. 
String $6$  starts
with bi-degree $(6,2)$ and ends in bi-degree $(2,6)$. That accounts for all the $2's$ in the image  of its path.  String $7$  starts
with bi-degree $(4,2)$ and ends in bi-degree $(2,4)$. 
\sas

Our computations using MAPLE yielded the following seven starters:
\sas

\noindent
$\DD_{1^5}$, $E_{2,0}\DD_{1^5}$,
$E_{3,0}\DD_{1^5}$, $E_{4,0}\DD_{1^5}$,
$3E_{2,0}E_{2,0}\DD_{1^5}+2EE_{3,0}\DD_{1^5}$,
$2E_{2,0}E_{2,0}E_{2,0}\DD_{1^5}-3EEE_{4,0}\DD_{1^5}$,
$E_{4,0}E_{2,0}\DD_{1^5}$.
\sas

Notice that for $n<7$ there is at most one starter at any bi-degree. If we find one that starts at bi-degree $(a,b)$ and ends at  $(b,a)$ we can safely complete the string in a construction of a basis. However, the result is not a basis that is consistent with Theorem I.1.
This is what happens with all the elements of the strings started by the $5^{th}$ and $6^{th}$ starters.

It will be good 
to make a few observations before we work on $n=6$.
Firstly, Remark 1.2 tells us that the operators $E_{r,0}$ commute whatever is the value of $r$. From Theorem I.1 we know that to obtain an alternant of bi-degree $(a,b)$ in $DH_n$ using  the polynomial 
$E_{r_1,0}E_{r_2,0}\cdots E_{r_b,0}\DD_{1^n}$ we
must require the equality 
$$
r_1+r_2+\cdots + r_b\ses 
\tttt{n \choose 2}-a
\eqno 2.2
$$
Due to the commutativity fact,  there is no loss in assuming that our sequences
are weakly decreasing 
($r_1\ge r_2\ge\cdots\ge r_b$) and satisfy 2.2. This is the number of $b$ parts partitions of $\tttt{n \choose 2}-a$.

Notice that our display at the row indexed by $t^5$ and column indexed by $q^2$ reveals that there is only one alternant in $DH_5$ with bi-degree $(5,2)$.
Yet the $2$ parts partitions of $5$ are $32$ and $41$.
Thus when we reduce the two polynomials 
$E_{3,0}E_{2,0}\DD_{1^5}$ and $E_{4,0}E_{1,0}\DD_{1^5}$ to have leading monomial with coefficient $1$ the resulting alternants must be identical!

This fact shows that we can only use partitions as an upper bound to the number of strings that start at a given bi-degree. However, in our computer experimentation with $DHA_6$ we discovered  that pairs that were discarded in $DHA_5$ had also to be discarded as factors in the construction of string starters in $DHA_6$. This suggests that the construction of a basis for $DHA_n$ might demand a recursion on $n$.

These  findings suggest that it might be possible to obtain a recursive construction of a basis.
The following result gives us a tool  for not using factors that have been discarded for $n-1$ in the construction of basis elements for $n$.
For  simplicity we will state it in the simplest useful case. 
\sas

\noindent{\bf Theorem 2.1}

{\ita Suppose that  $C_{n-1}(q,t)\big|_{t^c q^{b-1}}=k-1$ and $C_{n}(q,t)\big|_{t^a q^{b}}=k$, then the polynomial 

\vskip -.18in
$$
P(X_n;Y_n)\ses E_{r_1,0}E_{r_2,0}\cdots E_{r_b,0}\DD_{1^n} 
\eqno 2.3
$$

\vskip -.11in
\noindent
cannot be used as a starter in bi-degree $(a,b)$ if

\vskip -.2in
$$
a\ses n-1+c - r_b
\eqno 2.4
$$

\vskip -.15in
\noindent 
for any discarded pairs of  solutions of }

\vskip -.2in
$$
r_1+r_2+\cdots+r_{b-1} \ses 
\tttt{n-1 \choose 2}\sms c \,.
\eqno 2.5
$$

\vskip -.2in 
\noindent{\bf Proof}

\vskip -.03in
Since our conjecture requires that

\vskip -.2in
$$
r_1+r_2+\cdots+r_b\ses 1+2+\cdots +n-1 \sms a \,.
\eqno  2.6
$$

\vskip -.07in
\noindent
this can be rewritten in an  inductive way by relating for $n$  what we already obtained for $n-1$

\vskip -.2in
$$
0\ses r_1+r_2+\cdots+r_{b-1} \sms 
(\tttt{n-1 \choose 2}\sms c )\ses 
 n-1\sms a \sps c
\sms   r_b\,.
\eqno  2.7
$$

\vskip -.1in
\noindent
This proves 2.4.
 Our task is now to explore $n=6$

\vskip -.13in
\hfill $\includegraphics[height=1.6 in]{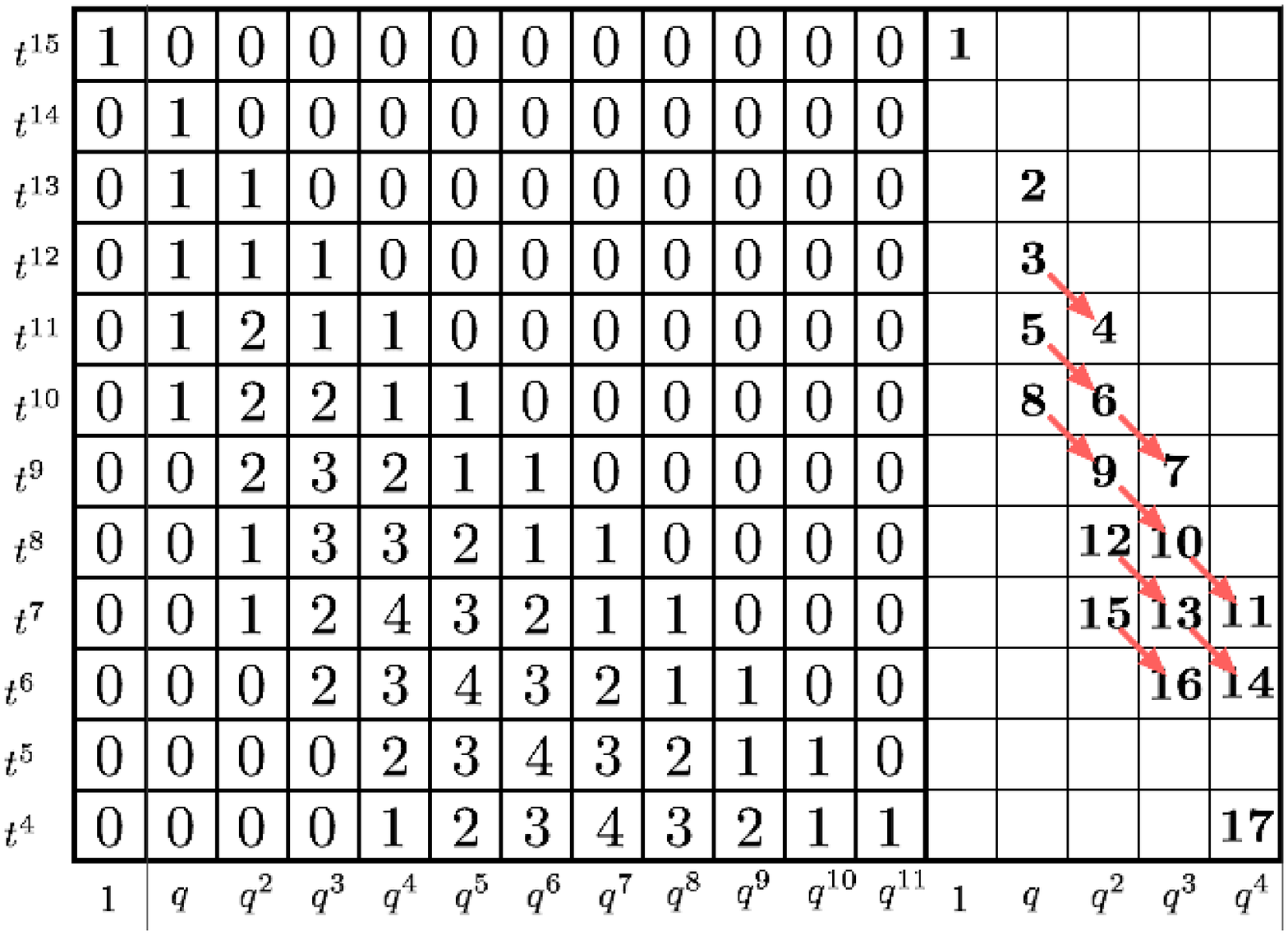}\ess\ess\ess$

\vskip -1.5in
\hsize=3.8in
 On the right we have a visual display of a portion of the Frobenius characteristic of the alternants in $DHA_6$. This is the polynomial $c_6(q,t)$, The number of strings start at bi-degree $(a,b)$ is again
\vskip -.18in
$$ 
c_6(q,t)\big |_{t^{a+1}q^{b-1}}
-\,
c_6(q,t)\big |_{t^{a }q^{b }}\le 1
\eqno 2.6
$$

\vskip -.06in
\noindent
Without the 0's the string starters 
are
$({\bf 1}\RA)\DD_{1^6}$, 
$({\bf 2}\RA)E_{2}\DD_{1^6}$,
\sas

\noindent
$({\bf 3}\RA)E_{3}\DD_{1^6}$, $({\bf 4}\RA)11 E_{2}E_{2}\DD_{1^6}+4 E_{3}E_{1}\DD_{1^6} $

\hsize=6.5in
\noindent  
$({\bf 5}\RA) E_{4}\DD_{1^6}$,  
$({\bf 6}\RA)2 E_{3}E_{2}\DD_{1^6}
+E_{4}E_{1}\DD_{1^6}$,
$({\bf 7}\RA)12E_{2}E_{2}E_{2}\DD_{1^6}+18E_{3}E_{2}E_{1}\DD_{1^6}+5E_{4}E_{1}E_{1}\DD_{1^6}
$
   
\sas

\noindent
$({\bf 8}\RA)E_{5}\DD_{1^6}$
$({\bf 9}\RA)3E_{3}E_{3}\DD_{1^6}+ 2E_{5}E_{1}\DD_{1^6}$, 
$({\bf10}\RA)28 E_{3}E_{2}E_{2}\DD_{1^6}-19E_{5}E_{1}E_{1}\DD_{1^6}$, 

\noindent
\bigsp\bigsp\bigsp $({\bf 11}\RA)5E_{2}E_{2}E_{2}E_{2}\DD_{1^6}+24E_{3}E_{2}E_{2}E_{1}\DD_{1^6}+
4E_{4}E_{2}E_{1}E_{1}\DD_{1^6}-8E_{5}E_{1}E_{1}E_{1}\DD_{1^6}$
\sas

\noindent
$({\bf 12}\RA)E_{4}E_{3}\DD_{1^6}$,
$({\bf 13}\RA)3E_{3}E_{3}E_{2}\DD_{1^6}+2E_{4}E_{3}E_{1}\DD_{1^6}$,

\bigsp\bigsp\bigsp\bigsp
$({\bf 14}\RA)20E_{3}E_{2}E_{2}E_{2}\DD_{1^6}+ 42E_{3}E_{2}E_{2}E_{1}\DD_{1^6}
+ 39E_{4}E_{2}E_{2}E_{1}\DD_{1^6}$
\sas

\noindent
$({\bf 15}\RA)E_{4}E_{4}\DD_{1^6}$,
$({\bf 16}\RA)5E_{3}E_{3}E_{3}\DD_{1^6}-9E_{4}E_{4}E_{1}\DD_{1^6}$
\sas

\noindent
$({\bf 17}\RA)E_{5}E_{2}E_{2}
\noindent
E_{2}\DD_{1^6}$.
\sas\sas

We can derive an algorithm for obtaining the number of starters for any $n$, that MAPLE permits, from the following calculation.
Since every string is of the form

\vskip -.15in
$$
P(X_n;Y_n)\RA EP(X_n;Y_n)\RA\cdots \RA E^kP(X_n;Y_n).
\eqno 2.7
$$

\vskip -.08in
\noindent
If $P(X_n;Y_n)$ is an alternant homogeneous 
of bi-degree $(u,v)$ then
the alternant $E^iP(X_n;Y_n)$ is
homogeneous of bi-degree $(u-i,v+i)$, with $E^kP(X_n;Y_n)$ homegeous of bi-degree $(v,u)$.
The contribution of this string to the Hilbert series of $DHA_n$ is the polynomial $\sum_{i=0}^{u-v} t^{u-i}q^{v+i}$. 
Note that 

\vskip -.16in
$$
t^u q^v\ses (qt)^{u+v\over 2}(\tttt{q\over t})^{v-u\over 2},
\eqno 2.8
$$

\vskip -.2in
\noindent
so we can write

\vskip -.2in
$$
\eqalign{
t^u q^v\Big(
1+&
(\tttt{q\over t})+\cdots + 
(\tttt{q\over t})^{u-v } \Big)
\ses
(qt)^{u+v\over 2}
(\tttt{t\over q})^{u-v\over 2}
\Big(1+ 
(\tttt{q\over t})+\cdots + 
(\tttt{q\over t})^{u-v } \Big)
\cr
&\ses
(qt)^{u+v\over 2}
(\tttt{t\over q})^{u-v\over 2}
{1 \sms (\tttt{q\over t})^{{u-v}+1}
\over
1\sms (\tttt{q\over t})
}
\ses (qt)^{u+v\over 2}
{(\tttt{t\over q})^{u-v\over 2} 
\sms (\tttt{q\over t})^{{u-v\over 2}+1}
\over
1\sms (\tttt{q\over t})
}
\cr
&\ses
(qt)^{u+v\over 2}
{(\tttt{t\over q})^{u-v \over 2}(\tttt{t\over q})^{1\over 2} 
\sms (\tttt{q\over t})^{{u-v\over 2}+1}(\tttt{t\over q})^{1\over 2}
\over
(\tttt{t\over q})^{1\over 2}\sms (\tttt{q\over t}) (\tttt{t\over q})^{1\over 2}
}
\ses
(qt)^{u+v\over 2}
{(\tttt{t\over q})^{u-v+1 \over 2} 
\sms (\tttt{q\over t})^{{u-v+1\over 2}}
\over
(\tttt{t\over q})^{1\over 2}\sms (\tttt{q\over t})^{1\over 2}
}
\cr}
\eqno 2.9
$$

\vskip -.1in
\noindent
Thus the Hilbert series of $DHA_n$ is the polynomial

\vskip -.16in
$$
h(q,t)\ses \sum_{u=0}^{n\choose 2}
 \sum_{v=0}^{n\choose 2}b_{u,v}
\,\,  (qt)^{u+v\over 2}
{(\tttt{t\over q})^{u-v+1 \over 2} 
\sms (\tttt{q\over t})^{{u-v+1\over 2}}
\over
(\tttt{t\over q})^{1\over 2}\sms (\tttt{q\over t})^{1\over 2}
}
\eqno 2.10
$$

\vskip -.1in
\noindent
where $b_{u,v}$ is the number of strings that start in bi-degree $(u,v)$. Making the specialization $t\RA q^{-1}$ gives

\vskip -.16in
$$
h(q,q^{-1})\ses \sum_{u=0}^{n\choose 2}
 \sum_{v=0}^{n\choose 2}b_{u,v}
\,\,  
{   q^{u-v+1  } 
\sms q^{{-(u-v+1) }}
\over
q\sms q^{-1}
}
$$

\vskip -.2in
\noindent
or better
$$
\ess\ess\ess (q\sms q^{-1})h(q,q^{-1})\ses \sum_{u=0}^{n\choose 2}
 \sum_{v=0}^{n\choose 2}b_{u,v}
\,\,  
 \big(  q^{u-v+1}   
\sms q^{-(u-v+1) }\big)
\eqno 2.11
$$

\noindent
We can rewrite the identity in 2.11 in the form
$$
\ess\ess\ess (q\sms q^{-1})h(q,q^{-1})\ses \sum_{u=0}^{n\choose 2}
 \sum_{v=0}^{n\choose 2}
\chi\big(u-v+1=r\big) 
 \big(  q^r   
- q^{-r }\big)\sum_{u-v+1=r}
b_{u,v}  
\eqno 2.12
$$
If we let $c_r=\sum_{u-v+1=r}
b_{u,v}$  then 2.12 becomes 
$$
\eqalign{  
\bu\bu\bu\bu(q\sms q^{-1})h(q,q^{-1})
 &\ses
 \sum_{u=0}^{n\choose 2}
 \sum_{v=0}^{n\choose 2}
\chi\big(u-v+1=r\big) 
 \,\big(  q^r   
- q^{-r }\big)
\,c_r \cr
&\ses
 \sum_{r=1}^{n\choose 2}
 \big(q^r-q^{-r }\big)
\, c_r
 \sum_{u=0}^{n\choose 2}
 \sum_{v=0}^{n\choose 2}
\chi\big(r=1-u+v \big)
\ses 
 \sum_{r=1}^{n\choose 2}
 \big(  q^r   
\sms q^{-r }\big)
\,c_r\,.
\cr} 
\eqno 2.13
$$

\vskip -.1in
Haiman has proved in [12] that $h(q,t)$ is exactly  the Garsia-Haiman $q,t$-Catalan. In [6] the specialization $t\RA q^{-1}$ is derived from Macdonald identities
to be related to the
none other than the q-analogue of the  number of Dyck paths in the $n\times n$ lattice rectangle, more precisely we have:
$$
q^{n \choose 2}h(q,q^{-1})\ses {1\over [n+1]_q}\Big[{2\,n \atop n}\Big]_q\, .
\eqno 2.14
$$
Thus the left hand side of 2.13 is
$$
 (q\sms q^{-1}) 
q^{-{n \choose 2}}
{q-1 \over q^{n+1}-1}
{(q^{n+1}-1)\cdots(q^{2\, n}-1)\over 
(q-1)\cdots(q^{ n}-1)}
\ses
q^{-{n \choose 2}-1}
{(q^{n+2}-1)\cdots(q^{2\, n}-1)\over 
(q^3-1)\cdots(q^{ n}-1)}
\eqno 2.15
$$
and 2.13  becomes 
$$
 q^{-{n \choose 2}-1}{(q^{n+2}-1)\cdots(q^{2\, n}-1)\over 
(q^3-1)\cdots(q^{ n}-1)}
\ses 
 \sum_{r=1}^{n\choose 2}
 \,\, c_r\,\, q^r   
\sms 
\sum_{r=1}^{n\choose 2}
c_r\,\, q^{-r }
\,
\eqno 2.16
$$
Since the number of starters is 
$\sum_{r=1}^{n\choose 2}c_r$ all we need is to sum the positive coefficients of the powers of $q$. Of course this is really only an algorithm, but we can also obtain a formula from 2.16.

To do this we first apply  an odd power of the Euler operator $q\del_q$ to both  side of 2.16$\,$ :
$$
\big(q\del_q\big)^{2k+1}\,\left ( q^{-{n \choose 2}-1}{(q^{n+2}-1)\cdots(q^{2\, n}-1)\over 
(q^3-1)\cdots(q^{ n}-1)}\right)
\ses 
 \sum_{r=1}^{n\choose 2}
 \,\,  r^{2k+1}\,c_r\, (q^{r}+q^{-r})  
\eqno 2.17
$$
The determinant of the matrix $\big\|z_r^{2k}\big\|_{r,k=1}^{{n(n-1)/2}+1}$ factorized  into a product of $(z_r^2-z_s^2)$ and each is different from zero as long as $z_r\neq z_s$.
this proves that the matrix
 $\big\|r^{2k+1}\big\|_{r,k=1}^{{n(n-1)/2}+1}$ has non zero determinant. Denoting by
 $\big\|d_{s,r}\big\|_{s,r=1}^{{n(n-1)/2}+1}$ the inverse,
it follows that
$$
\sum_{r=1}^{n\choose 2}d_{s,r}
\big(q\del_q\big)^{2k+1}\,\left ( q^{-{n \choose 2}-1}{(q^{n+2}-1)\cdots(q^{2\, n}-1)\over 
(q^3-1)\cdots(q^{ n}-1)}\right)
\ses 
 \sum_{r=1}^{n\choose 2}
 \,\,c_r\, (q^{r}+q^{-r})  
\eqno 2.18
$$
Thus
$$
\tttt{1\over 2}\sum_{r=1}^{n\choose 2}d_{s,r}
\big(q\del_q\big)^{2k+1}\,\left ( q^{-{n \choose 2}-1}{(q^{n+2}-1)\cdots(q^{2\, n}-1)\over 
(q^3-1)\cdots(q^{ n}-1)}\right)\Big|_{q=1}
\ses 
 \sum_{r=1}^{n\choose 2}
 \,\,c_r\,   
\eqno 2.19
$$

\vskip -.2in
\noindent
as desired.
\sas

Using the algorithm in 2.16, we obtain the following sequence giving the number of string starters for each $n\ge 2$. 
\sas

\noindent
$ \includegraphics[height=1.4 in]{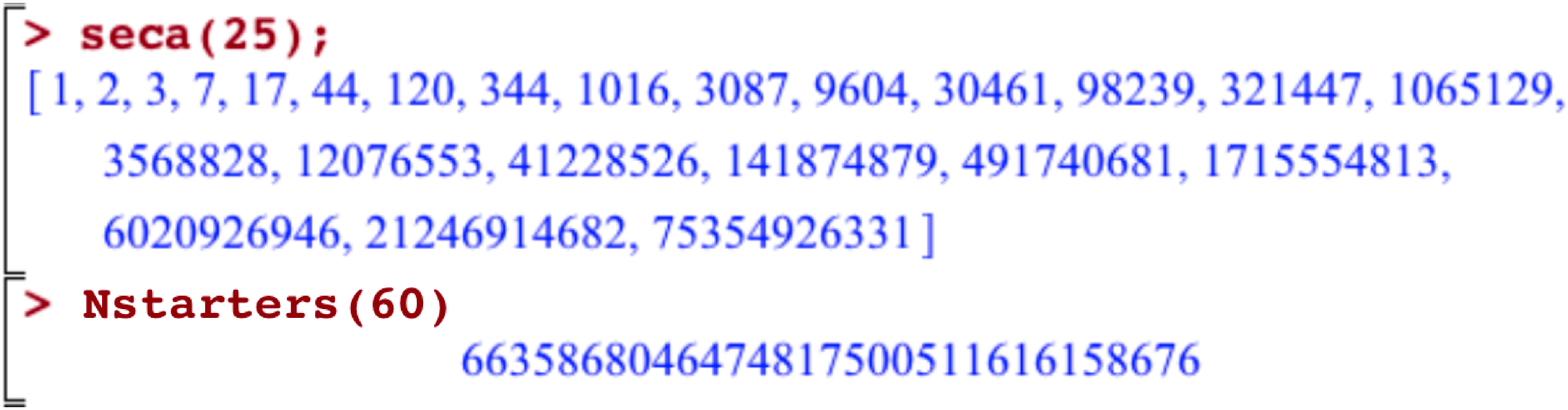}  $

\vskip -0.8 in

\hfill 2.20 

\vskip  0.5 in

\noindent
$\ess\ess\ess$The Encyclopedia of Integer Sequences was not aware of the existence of  this sequence.
\sas

Conjecture in [1] and Theorem 2.1, give us an algorithm for constructing starters, this done a basis for $DHA_n$ is easily constructed. From the display in 2.20, we see that the number of starters for $n=7$ is 
$44$ and for $n=8$ is already $120$.
With some work we may be able to construct the starters for
$n=7$, with even more work we may be able to do it for $n=8$. However, a paper [1] published in 1998 by Ed Allen  yields us a   basis for every $n$. To see how we can use [1]
we need only to depict an example.

 \hfill $ \includegraphics[height=1.3 in]{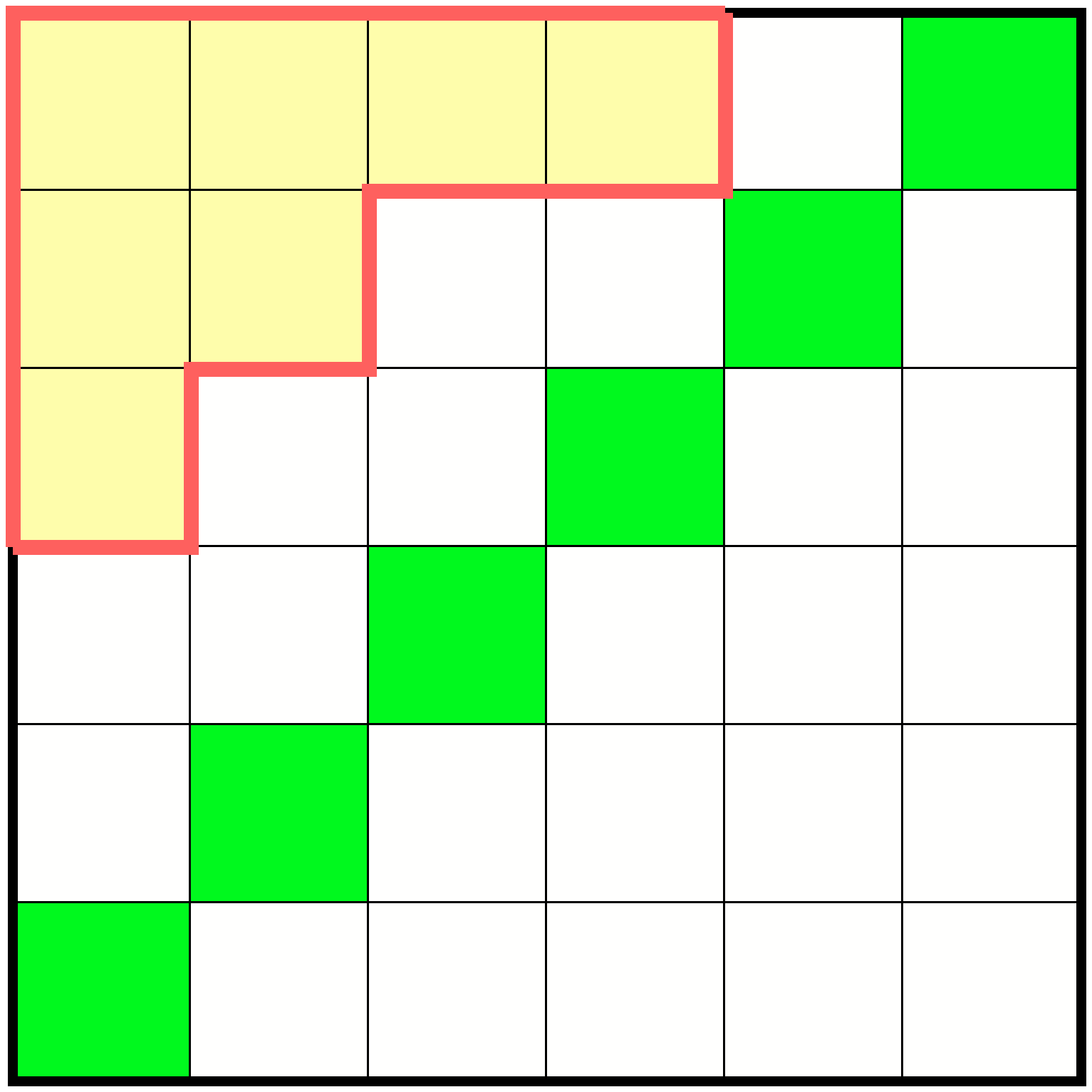}  $

\hsize 5in
\vskip -1.35in
On the right we have the lattice square $\CL_6$. The green cells are the  lattice diagonal. We have depicted in yellow the cells of the english Ferrers diagram of the partition $[4,2,1]$. Starting from $(0,0)$ we have $3$ black North steps, followed by $3$ red North steps, delimiting a Dyck path $D$ in $\CL_6$.
The $7$ yellow cells  give the co-area of $D$. We will call the yellow partition the ``{\ita co-partition}''  of $D$. In [6] we conjecture that the  dimension of $DHA_n$ is given by the number Dyck paths in $\CL_n$. Accordingly in [1] it is conjectured that the following succession of steps
constructs

\hsize 6.5in
\noindent
a basis of $DHA_n$. First step, compute the power function expansion of the  Schur function indexed by the co-partition for each Dyck path
in $\CL_n$. Second step, replace each $p_r$ by the operator $E_r$. Of course that means $p_r^k$ gets replaced by $E_r$ repeated  $k$ times. This done, apply the resulting differential operator to   the Vandermonde determinant in $\xon$. This sequence of steps should yield a  basis for $DHA_n$.

Since the power function expansion of a Schur function involves many terms
the resulting alternant may not be 
bi-homogeneous. For instance for the co-partition $[4,2,1]$ the power  basis expansion of the corresponding Schur function is the sum of the following terms
$$
\big[\tttt{1\over 144}p_1^7,
\tttt{-1\over 72}p_1^4p_3,
\tttt{-1\over 48}p_1^3p_2^2,
\tttt{-1\over 24}p_1^3p_4,
\tttt{-1\over 12}p_3p_1^2p_2,
\tttt{-1\over 18}p_1p_3^2,
\tttt{-1\over 24}p_3p_2^2,
\tttt{-1\over 24}p_3p_2^2,
\tttt{-1\over 24}p_3p_2^2,
\tttt{-1\over 24}p_3p_2^2,
\tttt{-1\over 24}p_3p_2^2,
\tttt{-1\over 24}p_3p_2^2\big]
\eqno 2.21
$$

The next result reinforces Conjecture in [1]. 
\sas

\noindent
{\bf Theorem 2.2}

{\ita On the validity of the conjecture in [1] it is possible to construct a bi-homogeneous basis for $DHA_n$.}

\noindent
{\bf Proof}

It suffices to  construct such a basis by an algorithm that works for all $n$. A fundamental fact is that $DHA_n$ is a bi-graded vector space. Moreover it is shown in [12] that $DHA_n$ has dimention equal to the number of Dyck paths in $\CL_n$. The conjecture in [1] starts with the power basis expansion of Schur functions indexed by co-partitions in $\CL_n$. In the next step we replace each $p_r$ in the resulting expansion by the operator $Er$ repeated as many times as the  exponent of $p_r$. This done we apply the 
resulting operator to the Vandermonde in $\xon$. The resulting polynomial will be a bi-homogeneous element of $DHA_n$, for every term in the power function expansion.

\vfill\supereject
\hfill $ \includegraphics[height=2 in]{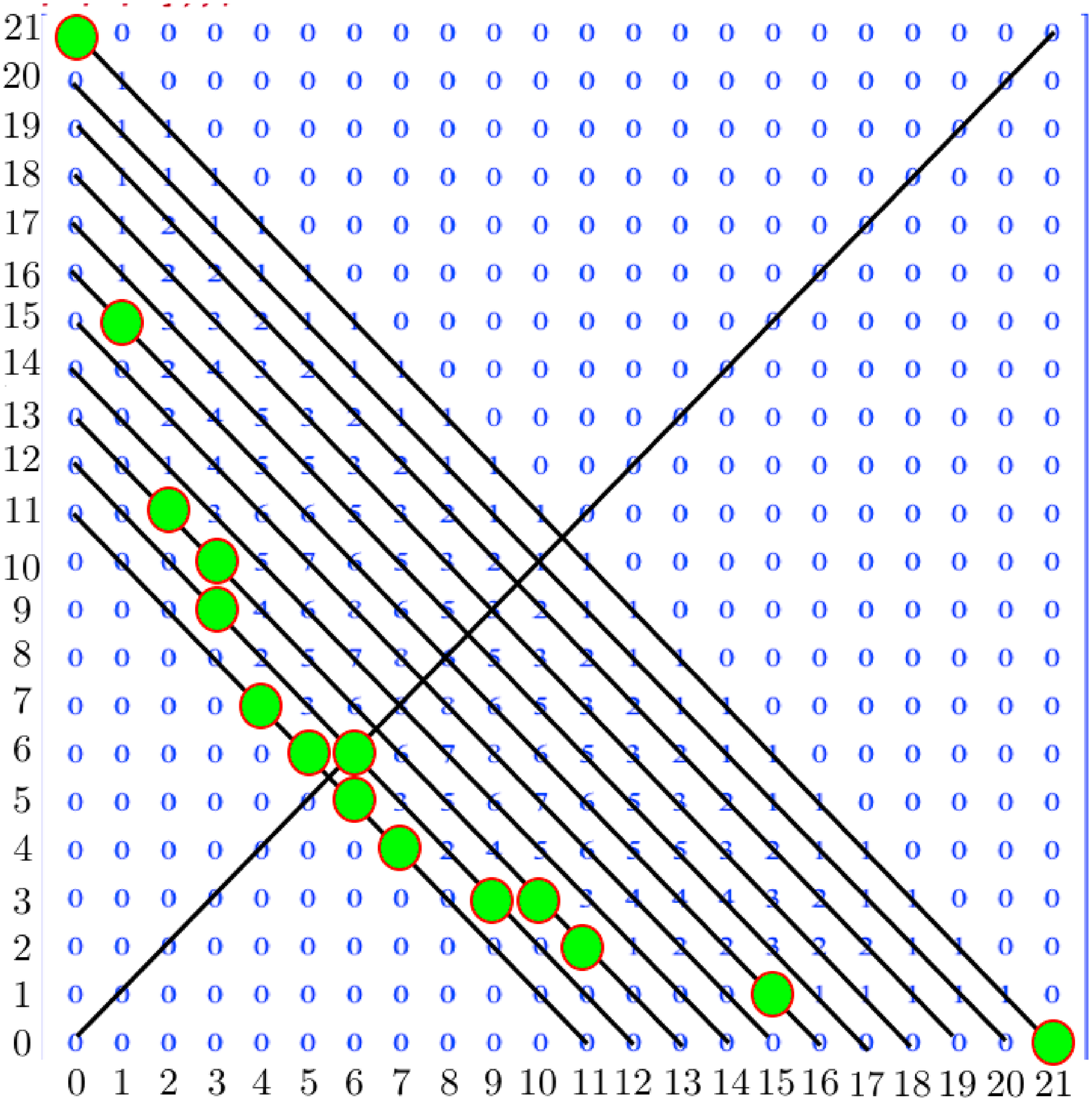}  $

\vskip -2.1 in
\hsize 4.4 in
To construct a basis we proceed with  $0\le i\le {n \choose 2}$,
starting by the diagonal of the string generated by the Vandermonde in $\xon$ in bi-degree ${n \choose 2},0$ and ending with the Vandermonde in $\yon$ in bi-degree $0,{n \choose 2}$.
Then process all the bi-degrees in the shorter diagonals as $i$ becomes progressively smaller.  

Say that in bi-degree $a,b$ the  Frobenius characteristic of $DHA_n$ tells
us that there are $m$ independent bi-homogeneous polynomials. To construct such a set we process from smallest to biggest the co-partition expansions breaking ties in the lex order of the partition sequence. We must eventually find $m$ independent polynomials of bi-degree $a,b$ if the basis conjectured in [1] is valid. This completes our proof.

\hsize 6.5 in

The display above  and on the right we exhibited  the Frobenius characteristic of $DHA_7$. The integers immediately to the left of the display give the  degrees in $x_i's$  of the  polynomials exhibited the rows. The integers at the bottom of the display give the  degrees in $y_i's$  of the  polynomials exhibited in the columns. The smaller integers
under the diagonals are the number of independent polynomials that occupy that particular bi-degree. the circles indicate the polynomials indexed by the partitions of $n=7$. For instance the circle in bi-degree $21,0$ is the Vandermonde in the $x_i's$ and the circle in bi-degree $0,21$ is the Vandermonde in the $y_i's$. Since our guide is the Frobenius characteristic of $DHA_n$, we can see that each of the four polynomials in  bi-degrees $(7,4), (6,5), (5,6), (4,7)$ will be picked up, in the construction of a basis, by our algorithm  on the validity of the conjecture in [1]  for $n=7$.
We mention this fact since N. Wallach (see [15]) that predicts, for any number $k$ of distinct variables, the polynomials that will occur in any basis. What we are witnessing here, for $k=2$ and $n=7$, the validity of a general result of this general result.
\sas

\noindent
{\bf Theorem 2.3}

{\ita In $DHA_n$ there is a basis of $ker \,F$ of the form
$$
E_2^{a_2}E_3^{a_3}\cdots E_{n-1}^{a_{n-1}}\DD_{1^n}\, +\, E_1\CL_{\bf Q}[E_1,E_2,\cdots ,E_{n-1}]\DD_{1^n} 
\eqno 2.22
$$
where $a_i\in Z_{\ge 0}$ and $\CL_{\bf Q}$ denotes the linear span with coefficients in $\bf Q$.}

\noindent{ \bf Proof}
\sas

It is well known from $sl[2]$ theory that If $W$ is a subspace of ${\bf C}[\xon;\yon]$ that is invariant under $E$ and $F$ then
$$  
W\ses E_1W\, \oplus\, \ker \, F|_{W}.  
\eqno 2.23
$$
Moreover any irreduble representation of dimension
$k+1$ starts with a polynomial $v_0$ such that
$$
 a)\ess\ess F\, v_0=0,
\bigsp\bigsp
 b)\ess\ess H\,v_0\ses -k\,  v_0.
\eqno 2.24
$$
Then  a basis of the representation is
$$
v_0,E v_0,E^2 v_0,\cdots,E^kv_0
\eqno 2.25
$$
Thus if $f\in DHA_n$ then
$$
f\ses\CL_{\bf Q}[E_1,E_2,\cdots,E_{n-1}]\DD_{1^n}
\, +\, E_1\CL_{\bf Q}[E_1,E_2,\cdots ,E_{n-1}]\DD_{1^n}.
\eqno 2.26
$$
From the  operator theorem it follows that there is a basis of the form
$$
\CL_{\bf Q}[E_2,E_3,\cdots,E_{n-1}]\DD_{1^n}.
\eqno 2.27
$$
If $F\, f=0$  then  we can assure that
$$
f\ses \CL_{\bf Q}[E_2,E_3,\cdots,E_{n-1}]\DD_{1^n} \,+\,  E_1\CL_{\bf Q}[E_1,E_2,\cdots ,E_{n-1}]\DD_{1^n}.
\eqno 2.28
$$
This proves the Theorem.

\page
\noindent
Acknowledgement:  The authors would like to thank Nolan Wallach for
his contributions to the mathematics in this paper.  The authors would
also like to thank Marino Romero for numerous conversations about the
contents of the paper.
\vskip .2in

\centerline{\bol Bibliography}
\sas
\item {[1]} E. Allen, {\ita A conjecture of a basis for the diagonal harmonic alternants
}, Discrete Mathematics {\bf 193} (1998) 33-42.
\ssas

\item {[2]} E. Carlsson and A. Oblomkov, {\ita Affine Schubert calculus and double coinvariants
}, arXiv 1801.09033v3, (2019).
\ssas

\item  {[3]} O, Egecioglu,  A. M. Garsia, {\bf Lectures in Algebraic  Combinatorics},
Springer Lecture Notes in Mathematics {\bf  2277} (2020)
\ssas

\item  {[4]} A. M. Garsia, G. Xin and M. Zabrocki, {\ita  Hall-Littlewood operators in the Theory of Parking Functions and Diagonal Harmonics},
International  Mathematical Research Notices, {\bf 2012 (6)} (2012), 1264--1299.
\ssas

\item {[5]} A. M. Garsia and M. Haiman,
{\ita Some Natural bigraded $S_n$-Modules and $q,t$-Kostka Coefficients}, 
Electronic J. of Combinatorics, Volume 3, Issue 2 (199) (The Foata Festschrift Volume). 
\ssas

\item {[6]} A. M. Garsia and M. Haiman,
{\ita A Remarkable  $q,t$-Catalan sequence and  $q$-Lagrange Inversion},
 J. Algebraic Combin., {\bf 5 (3)} (1996), 191--244.
\ssas

\item {[7]} A. M. Garsia and M. Haiman and G.Tesler,
{\ita Explicit Plethystic Formulas for the Macdonald q,t-Kostka  Coefficients}, {\bf Seminaire Lotharingien de Combinatoire, B42m, 42 pp. (1999).} 
\ssas

\item {[8]} A. M. Garsia and J. Haglund,
{\ita A positivity result in the theory of Macdonald polynomials}, Communicated by Ronald L. Graham, University of California San Diego, January26, 2001. PNAS April 10 98 (8). (2001)  4313-4316.
\ssas

\item {[9]} A. M. Garsia and J. Haglund,
{\ita A proof of the $q,t$-Catalan positivity conjecture},
Discrete Mathematics, {\bf 256 (2)} (2002), 677--717.
\ssas

\item {[10]} J. Haglund, M. Haiman, N. Loehr, J. B. Remmel and A. Ulyanov,
{\ita A combinatorial formula for the character of the diagonal coinvariants},
{ Duke Math. J.}, {\bf 126} (2005), 195--232.
\ssas

\item  {[11]}  J. Haglund, {\ita The q,t-Catalan numbers and the space of diagonal harmonics}, volume 41 of
University Lecture Series. American Mathematical Society, Providence, RI (2008).
\ssas

\item  {[12]}
M. Haiman.
{\ita Hilbert schemes, polygraphs and the Macdonald positivity conjecture.}
J.  Amer. Math. Soc. {\bf 14} (2001), 941-1006.
\ssas

\item  {[13]}
M. Haiman.
{\ita Vanishing theorems and character formulas for the  Hilbert scheme of points in the plane.}
Invent. Math. {\bf 149} (2002), 371--407.
\ssas 

\item {[14]}  I. G. Macdonald,
{\ita Symmetric functions and Hall polynomials},
 2nd Ed. Reprint of the 2008 paperback edition,
 Oxford University Press, New York (2015).
\ssas 

\item {[15]} Nolan Wallach, {\ita The representation of
$GL(k)$ on the alternants of minimal degree for the diagonal action of $S_n$ on $k$ copies of the permutation representation
,}  math arXiv September 6 (2021).
\ssas

\item {[16]} Hermann Weyl, {\bf
The Classical Groups Their Invariants and Representations}, Princeton Publications in Mathematics.

\end